\documentclass[3p,10pt,a4paper,twoside,fleqn,sort&compress]{filomat}
\usepackage{amssymb,amsmath,latexsym}
\usepackage[varg]{pxfonts}

\newtheorem{theorem}{Theorem}[section]

\newtheorem{lemma}[theorem]{Lemma}

\newtheorem{remark}[theorem]{Remark}

\newcommand{\FilVol}{xx}
\newcommand{\FilYear}{yyyy}
\newcommand{\FilPages}{zzz--zzz}
\newcommand{\FilDOI}{(will be added later)}
\newcommand{\FilReceived}{Received: dd Month yyyy; Accepted: dd Month yyyy}

\begin{document}

\title{Statistical approximation properties of Stancu type $q$-Baskakov-Kantorovich operators}

\author[affil1,affil2]{Vishnu Narayan Mishra}
\ead{vishnunarayanmishra@gmail.com}
\author[affil1]{Preeti Sharma}
\ead{preeti.iitan@gmail.com}

\author[affil3]{Adem Kili\c{c}man}
\ead{akilicman@putra.upm.edu.my}

\author[affil1,affil4]{Dilip Jain}
\ead{dilip18pri@gmail.com}

\address[affil1]{Department of Applied Mathematics \& Humanities,
Sardar Vallabhbhai National Institute of Technology, Ichchhanath Mahadev Dumas Road, Surat -395 007 (Gujarat), India}

\address[affil2]{L. 1627 Awadh Puri Colony Beniganj, Phase -III, Opposite - Industrial Training Institute (I.T.I.), Ayodhya Main Road \newline Faizabad-224 001, (Uttar Pradesh), India}

\address[affil3]{Department of Mathematics, Universiti Putra Malaysia (UPM), 43400 Serdang, Selangor, Malaysia}

\address[affil4]{Ward No. 12, Sumerganjmandi, Indergarh, (Bundi), Rajasthan - 323614, India }
\newcommand{\AuthorNames}{Mishra et al.}

\newcommand{\FilMSC}{Primary $41A10$, $41A25$; Secondary $41A36$}
\newcommand{\FilKeywords}{($q$-integers, $q$-Baskakov-Kantorovich operators, rate of statistical convergence, modulus of continuity, Lipschitz type functions)}
\newcommand{\FilCommunicated}{ Prof. H.M. Srivastava}
\newcommand{\FilSupport}{Ministry of Human Resource Development, India}

\begin{abstract}
In the present paper, we consider  Stancu type generalization of Baskakov-Kantorovich operators based on the $q$-integers and obtain statistical and weighted statistical approximation properties of these operators. Rates of statistical convergence by means of the modulus of continuity and the Lipschitz type function are also established for said operators. Finally, we construct a bivariate generalization of the operator and also obtain the statistical approximation properties. 
\end{abstract}

\maketitle

\makeatletter
\renewcommand\@makefnmark%
{\mbox{\textsuperscript{\normalfont\@thefnmark)}}}
\makeatother

\section{Introduction}
In the last decade, some new generalizations of well known positive linear operators based on $q$-integers were introduced and studied by several authors.~Our aim is to investigate statistical approximation properties of a Stancu type $q$-Baskakov-Kantorovich operators.~Firstly, Baskakov-Kantorovich operators based on $q$-integers was introduced by Gupta and Radu in \cite{VR} and they established some approximation results.\\
\indent Later, I. B\"{u}y\"{u}kyazici  and Atakut \cite{AB} introduced a new Stancu type generalization of $q$-Baskakov operators which is defined as
\begin{equation}\label{a}
{{\mathcal{L}^{(\alpha , \beta)}_n}}(f; q, x) = \sum_{k = 0}^{\infty} q^{\frac{k(k-1)}{2}} \frac{{D^k_q} (\phi_n (x))}{[k]_q ! } (-x)^k f\bigg( \frac{1}{q^{k-1}} \frac{[k]_q + q ^{k-1} \alpha}{[n]_q + \beta } \bigg),
\end{equation}
where $0 \leq \alpha \leq \beta ,~ q \in (0, 1), ~f \in C[0, \infty)$.\\ 
\indent Let $\{ \phi_n \} ~  (n = 1,2,3, \ldots ),~ \phi_n : \mathbb{R} \rightarrow \mathbb{R}$ be sequence which satisfies following conditions:
\begin{itemize}
\item[(i)]  $\phi_n ~(n = 1,2,3, \ldots )$ $k$-times continuously $q$-differentiable in any closed interval $[0, A],$ where $ A> 0,$
\item[(ii)] $\phi_n(0) = 1,~(n = 1,2, \ldots ),$
\item[(iii)] for all $x \in [0, A]$ and $(k = 1,2,\ldots; ~ n = 1,2, \ldots ),~ (-1)^k {D^k_q} (\phi_n(x)) \geq 0, $
\item[(iv)] there exist a positive integer $m(n)$, such that 
\begin{equation*}
{D^k_q}(\phi_n (x)) = - [n]_q {D^{k-1}_q} \phi_{m(n)}(x),~ (k = 1,2,\ldots;~ n = 1,2, \ldots ),
\end{equation*}
\item[(v)] $\lim\limits_{n \rightarrow \infty} \frac{[n]_q}{[m(n)]_q} = 1.$ 
\end{itemize}
\indent We first start by recalling some basic definitions and notations of $q$-calculus.\,We consider $q$ as a real number satisfying $0<q<1$.\\
For each non negative integer $n$, we define the $q$-integer $[n]_q$ as
\begin{equation*}
\displaystyle [n]_q = \left\{
\begin{array}{ll}
\frac{1-q^n}{1-q}, & \hbox{$q\neq1$}, \\
n,& \hbox{$q=1$}.
\end{array} \right.
\end{equation*}
The $q$-factorial is defined as
\begin{equation*}
[n]_q!=\left\{
\begin{array}{ll}
[n]_q[n-1]_q[n-2]_q...[1]_q,  & \hbox{$n=1,2,...$}, \\
1,& \hbox{$n=0$.}
\end{array}\right. \end{equation*}
We observe that
\begin{equation*} 
(1+x)_q^n=(-x;q)_n= \left\{ \begin{array}{ll} (1+x)(1+qx)(1+q^2x)...(1+q^{n-1}x), & \hbox{$n=1,2,...$},\\
1,& \hbox{$n=0.$} \end{array}\right.
\end{equation*} 
Also, for any real number $\alpha$, we have
$$(1+x)_q^\alpha =\frac{(1+x)_q^\infty}{(1+q^\alpha x)_q^\infty}.$$ 
In special case, when $\alpha$ is a whole number, this definition coincides with the above definition.\\
The $q$-binomial coefficients are given by 
$$\binom{n}{k}_q = \frac{[n]_q!}{[k]_q![n-k]_q!}, \,~ 0\leq k \leq n.$$
\newline
The $q$-derivative $D_q f$ of a function $f$ is given by
\begin{equation*} D_q(f(x))=\frac{f(x)-f(qx)}{(1-q)x},~ x\neq 0.
\end{equation*}
The $q$-Jackson integral is defined as
\begin{equation*}
\int_0^a f(x) d_qx=(1-q)a\sum_{n=0}^\infty f(aq^n)q^n, ~a>0.
\end{equation*}
Over a general interval $[a,b],~ 0 < a < b,$ one defines

\begin{equation*}
\int_a^b f(x) d_qx=\int_0^b f(x) d_qx - \int_0^a f(x) d_qx.
\end{equation*}
Throughout the paper, we use $e_i$ the test functions defined by $e_i(t) := t^i,$  where $i=0,1,2.$
First we need the following auxiliary result.\\
\newline
\indent Let $ \{ { \phi_n } \} $
be a sequence of real functions on $ \mathbb{R}_+ = [0, \infty )$ which are $k$-times continuously $q$-differentiable on $\mathbb{R}_{+}$ satisfying following conditions:
\begin{itemize}
\item[($p_1$)] $\phi_n(0) = 1,~~(n = 1,2, \ldots ),$
\item[($p_2$)] for $k \in \mathbb{N}_0 =\mathbb{N}\cup \{0 \}$ and $n \in \mathbb{N},~ (-1)^k {D^k_q} (\phi_n(x)) \geq 0,~ x \in \mathbb{R_+},$
\item[($p_3$)] there exist a positive integer $m(n)$, such that 
\begin{equation*}
{D^k_q}(\phi_n (x)) = - [n]_q {D^{k-1}_q} \phi_{m(n)}(x), ~(k = 1,2,\ldots;~ n = 1,2, \ldots ),
\end{equation*}
\item[($p_4$)] $\lim\limits_{n \rightarrow \infty} \frac{[n]_q}{[m(n)]_q} = 1.$
\end{itemize}

 ~\indent Under the condition ($p_1$) - ($p_4$),  \c{C}. Atakut and   \.{I}. B\"{u}y\"{u}kyazici \cite{AI} defined a new generalization of Stancu type $q$-Baskakov-Kantorovich operators as follows
 
\begin{equation}\label{b}
{\mathcal{L}^{*(\alpha, \beta)}_n}(f;q,x) = ([n]_q + \beta) \sum_{k = 0}^{\infty} q^{k(k-1)/2} \frac{{D^k_q} (\phi_n (x))}{[k]_q ! } (-x)^k \int_{q \left( \frac{[k]_q + q^{k-1}\alpha}{[n]_q+ \beta}\right)}^{\frac{[k+1]_q + q^k \alpha}{[n]_q+ \beta}}{ } f(q^{-k+1} t) d_qt, 
\end{equation}
where $x \in \mathbb{R_+},~n \in \mathbb{N}, ~ 0 \leq \alpha \leq \beta.$\\
\indent To obtain the approximation results we need the following Lemmas in what follows.

\begin{lemma}\cite{AB}
${\mathcal{L}}^{(\alpha, \beta)}_n$ be defined by (\ref{a}). Then the following identities hold
\begin{align*}
{\mathcal{L}}^{(\alpha , \beta)}_n(e_0; q, x) &= 1,\\ 
{\mathcal{L}}^{(\alpha , \beta)}_n(e_1; q, x) &= \frac{[n]_q}{[n]_q + {\beta}}x + \frac{\alpha}{[n]_q+ \beta},\\ 
{\mathcal{L}}^{(\alpha , \beta)}_n(e_2; q, x) &= \frac{[n]_q [m(n)]_q}{q([n]_q + {\beta})^2}x^2 + \frac{[n]_q (2 \alpha +1)}{([n]_q + {\beta})^2}x + \frac{{\alpha}^2}{({[n]_q + {\beta}})^2}.
\end{align*}
\end{lemma}

\begin{lemma}\cite{AI} The following relations are satisfied:
\begin{align*}
\int_{q \left( \frac{ [k]_q + q^{k-1}\alpha}{[n]_q+ \beta}\right)}^{\frac{[k+1]_q + q^k \alpha}{[n]_q+ \beta}} d_qt
& =  \frac{1}{[n]_q+ \beta},\\
 \int_{q \left( \frac{ [k]_q + q^{k-1}\alpha}{[n]_q+ \beta}\right)}^{\frac{[k+1]_q + q^k \alpha}{[n]_q+ \beta}}td_qt &= \frac{[2]_q[k]_q + q^k (1+ 2 \alpha)}{[2]_q{([n]_q+ \beta)}^2},\\
\int_{q \left( \frac{ [k]_q + q^{k-1}\alpha}{[n]_q+ \beta}\right)}^{\frac{[k+1]_q + q^k \alpha}{[n]_q+ \beta}} t^2d_qt &=\frac{[3]_q{[k]_q^2} + q^k [k]_q\big((1+ 3 \alpha) [2]_q + 1 \big) + (1 + 3 \alpha + 3 {\alpha}^2) q^{2k}}{[3]_q ([n]_q+ \beta)^3}.
\end{align*}
\end{lemma}

\begin{lemma}\label{c}\cite{AI}
Let $e_{i}= t^i$, where $i=0,1,2.$  For all $x\in\mathbb{R_+},$ $n\in \mathbb{N},$  $\alpha, \, \beta\geq 0$
and $0<q<1,$ we have
\begin{align*}
{\mathcal{L}^{*(\alpha, \beta)}_n}(e_0;q,x) &= 1,\\
{\mathcal{L}^{*(\alpha, \beta)}_n}(e_1;q,x)&= \frac{[n]_q}{[n]_q + \beta} x + \frac{q(1 + 2 \alpha)}{[2]_q([n]_q + \beta)},\\
{\mathcal{L}^{*(\alpha, \beta)}_n}(e_2;q,x)&= \frac{[n]_q [m(n)]_q}{q([n]_q + {\beta})^2}x^2 + \frac{[n]_q \big[[3]_q+q\big((1+3\alpha)[2]_q+1\big)\big]}{[3]_q([n]_q + {\beta})^2}x + \frac{q^2(1+ 3 \alpha + 3{\alpha}^2)}{[3]_q({[n]_q + {\beta}})^2}.
\end{align*}
 \end{lemma}

\begin{remark}
From Lemma (\ref{c}), we have
\begin{align*}
{\alpha}_n(x)={\mathcal{L}^{*(\alpha, \beta)}_n}(t-x;q,x)&= \left(\frac{[n]_q}{[n]_q + \beta}-1\right)x+\frac{q_n(1+2\alpha)}{[2]_q([n]_q+\beta)},\\ 
{\delta}_n(x)= {\mathcal{L}^{*(\alpha, \beta)}_n}((t-x)^2;q,x)& = \left (\frac{[n]_q [m(n)]_q}{q({[n]_q + \beta)}^2}+ 1- \frac{2[n]_q}{[n]_q + \beta} \right) x^2 \\&~~+ \bigg(\frac{[n]_q \big[ [3]_q+ q\big((1+3\alpha)[2]_q+1\big)\big] }{{ ([n]_q + \beta) }^2[3]_q} - \frac{2q(1+2\alpha)}{[2]_q([n]_q + \beta)} \bigg) x+ \bigg( \frac{q^2(1+ 3 \alpha+ 3 {\alpha}^2)}{[3]_q {([n]_q + \beta)}^2 } \bigg).
\end{align*}
\end{remark}

\begin{remark}
If we put $q$=1, we get the moment of Stancu type Baskakov-Kantorovich operators as
\begin{align*}
&{\mathcal{L}^{*(\alpha, \beta)}_n}(e_1;1,x)= \frac{n}{(n + \beta)} x + \frac{(1 + 2 \alpha)}{2(n + \beta)},\\
&{\mathcal{L}^{*(\alpha, \beta)}_n}(e_2;1,x)= \frac{n\, m(n)}{(n + {\beta})^2}x^2 + \frac{ 2n( \alpha +1)}{(n + {\beta})^2}x + \frac{1+ 3 \alpha + 3{\alpha}^2}{3({n + {\beta}})^2},\\
&{\mathcal{L}^{*(\alpha, \beta)}_n}(t-x;1,x)= \left(\frac{n}{(n + \beta)}-1\right)x+\frac{(1+2\alpha)}{2(n+\beta)},\\
&{\mathcal{L}^{*(\alpha, \beta)}_n}((t-x)^2;1,x)= \left (\frac{n\,m(n)}{({n+ \beta)}^2}+ 1- \frac{2 n}{(n + \beta)} \right) x^2 + \left(\frac{2n (1+ \alpha) }{{ (n + \beta) }^2 } - \frac{(1+2\alpha)}{(n+ \beta)} \right) x + \left( \frac{(1+ 3 \alpha+ 3 {\alpha}^2)}{3 {(n + \beta)}^2 } \right).
\end{align*}
\end{remark}

\section{Korovkin type statistical approximation properties}
 
The idea of statistical convergence was introduced independently by Steinhaus \cite{ST}, Fast \cite{HF} and Schoenberg \cite{98}. The study of the statistical convergence for sequences of linear positive operators was attempted in the year 2002 by A.D. Gadjiev and C. Orhan \cite{GO}.
Recently the idea of statistical convergence has been used in proving some
approximation theorems. It was shown that the statistical versions are stronger than the classical ones.
Authors have used many types of classical operators and test functions to
study the Korovkin type approximation theorems which further motivate to
continue the study. In particular, Korovkin type approximation theorems
\cite{krkn} was proved by using statistical convergence by various authors, e.g. \cite{HM3,HM1,SG1,HK,ms214,FA,HM}. In the recent years, Stancu type generalization of the certain operators and trigonometric approximation of signals by different types of summability operators have been studied by several other researchers, we refer some of the important papers in this direction as (\cite{M1}-\cite{BA}) etc.\\
Korovkin type approximation theory has also many useful connections, other than classical approximation theory, in other branches of mathematics (see Altomare and Campiti in \cite{ac22}).\\
~\indent Now, we recall the concept of statistical convergence for
sequences of real numbers which was introduced by Fast \cite{HF} and
further studied by many others.

\parindent=8mm Let $K\subseteq \mathbb{N}$ and $K_{n}=\left\{ j\leq n:j\in
K\right\}.$ Then the $natural~density$ of $K$ is defined by $\delta
(K)={\lim\limits_{n}}~ n^{-1}|K_{n}|$ if the limit exists, where $|K_{n}|$ denotes the
cardinality of the set $K_{n}$.

\parindent=8mm A sequence $x=(x_{j})_{j\geq1}$ of real numbers is said to be $%
statistically$ $convergent$ to $L$ provided that for every $\epsilon >0$ the
set $\{j\in \mathbb{N}:|x_{j}-L|\geq \epsilon \}$ has natural density zero,
i.e. for each $\epsilon >0$,
\begin{equation*}
\lim\limits_{n}\frac{1}{n}|\{j\leq n:|x_{j}-L|\geq \epsilon \}|=0.
\end{equation*}
It is denoted by $st-\lim\limits_{n}x_{n}=L$.\\
We consider a sequence $q=(q_{n}),$ $q_{n}\in $ $(0,1),$ such that
\begin{equation}\label{a1}
\lim\limits_{n\rightarrow \infty }q_{n}=1.
\end{equation}
\newline
The condition (\ref{a1}) guarantees that $[n]_{q_{n}}\rightarrow
\infty $ as $n\rightarrow \infty .$
\newline
Now, let us recall the following theorem given by Gadjiev and Orhan \cite{GO}.
\begin{theorem}
If the sequence of linear positive operators $A_n:C_{M}[a,b]\rightarrow C[a,b]$ satisfies the conditions
\begin{equation}\label{wa1}
st-\lim_{n}\|A_n(e_\nu;\cdot)- e_\nu\|_{C[a,b]}=0,~~e_{\nu}(t)= t^{\nu} \text{ for }  \nu=0,1,2,
\end{equation}
then, for any function $f \in C_{M}[a,b]$, we have
$$st-\lim\limits_{n} \| A_n(f;\cdot)- f \|_{C[a,b]}=0,$$
where $C_M[a,b]$ denotes the space of all functions $f$ which are continuous in $[a,b]$ and bounded on the all positive axis.
\end{theorem}
In \cite{DO} Do\u{g}ru and Kanat defined the Kantorovich-type modification of Lupa\c{s} operators as follows:
\begin{equation}\label{dk}
\tilde{R}_n(f;q;x)=[n+1]\sum_{k=0}^{n}\bigg(\int_{\frac{[k]}{[n+1]}}^{\frac{[k+1]}{[n+1]}} ~ f(t) d_{q}t\bigg)\left(\begin{array}{c}n \\k \end{array}\right) \frac{q^{-k}q^{k(k-1)/2} x^k (1-x)^{(n-k)}}{(1-x+qx)\cdots(1-x+q^{n-1} x)}.
\end{equation}
Do\u{g}ru and Kanat \cite{DO} proved the following statistical Korovkin-type approximation
theorem for operators (\ref{dk}).

\begin{theorem}
Let $q:=(q_n),~ 0 < q < 1$, be a sequence satisfying the following conditions:\\
\begin{equation}\label{d}
 st-\lim_n q_n = 1,~ st - \lim_n q_n^n = a ~( a < 1)~ and~ st - \lim_n \frac{1}{[n]_q} = 0,
\end{equation}
then if $f$ is any monotone increasing function defined on $[0, 1]$, for the positive linear operator $\tilde{R}_n(f;q;x)$, then
$$ st- \lim_n {\Vert \tilde{R}_n(f;q;\cdot) - f \Vert}_{C[0, 1]} = 0$$ holds. 
\end{theorem}
In \cite{do26} Do\u{g}ru gave some examples so that $(q_{n})$ is statistically convergent to $1$ but it may not convergent to $1$ in the ordinary case.


\begin{theorem} \label{A1}
Let ${{\mathcal{L}^{*(\alpha , \beta)}_n}}$ be the sequence of the operators (\ref{b}) and the sequence
$q = (q_n)$ satisfies (\ref{d}). Then for any function $f \in C[0, A] \subset C[0,\infty),~ A > 0 ,$ we have
\begin{equation}
st- \lim_n \Vert {{\mathcal{L}^{*(\alpha , \beta)}_n}}( f; q, \cdot) - f \Vert = 0,
\end{equation}
where $C[0,A]$ denotes the space of all real bounded functions $f$ which are continuous in
$[0,A].$
\end{theorem}
\textbf{Proof.}
Let $e_{i}= t^i,$ where $i=0,1,2.$   Using ${\mathcal{L}^{*(\alpha, \beta)}_n}(1;q_n,x)=1,$ it is clear that\\
$$st-\lim_{n}\|{\mathcal{L}^{*(\alpha, \beta)}_n}(1;q_n,x)-1\|=0.$$
Now by Lemma (\ref{c})(ii), we have
\begin{equation*}
\lim_{n\rightarrow \infty }\Vert L_{n}^{\ast(\alpha,\beta)}(t;q_{n},x)-x\Vert
=\left\|\frac{[n]_q}{[n]_q + \beta} x + \frac{q(1 + 2 \alpha)}{[2]_q([n]_q + \beta)}-x\right\|\leq \frac{\beta}{[n]_q + \beta} x + \frac{q(1 + 2 \alpha)}{[2]_q([n]_q + \beta)}.
\end{equation*}
 For given $\epsilon >0$, we define the following sets:
 \begin{equation*}
 S=\{k:\Vert{\mathcal{L}^{*(\alpha, \beta)}_n}(t;q_{k},x)-x\Vert\geq \epsilon\},
 \end{equation*}
and 
\begin{equation}\label{e}
S'=\left\{k: \frac{\beta}{[k]_q + \beta} x + \frac{q(1 + 2 \alpha)}{[2]_q([k]_q + \beta)}\geq \epsilon\right\}.
\end{equation}
It is obvious that $S\subset S^{\prime }$, it can be written as
\begin{equation*}
\delta \left( \{k\leq n:\Vert L_{n}^{\ast(\alpha,\beta)}(t;q_{k},x)-x\Vert \geq
\epsilon \}\right) \leq \delta \left( \{k \leq n: ~ \frac{\beta}{[k]_q + \beta} x + \frac{q(1 + 2 \alpha)}{[2]_q([k]_q + \beta)}\geq \epsilon \}\right).
\end{equation*}
By using (\ref{d}), we get
$$st-\lim_{n}\left(\frac{\beta}{[n]_q + \beta} x + \frac{q(1 + 2 \alpha)}{[2]_q([n]_q + \beta)}\right)=0.$$ 
So, we have
$$\delta \bigg(\big\{k \leq n: ~ \frac{\beta}{[k]_q + \beta} x + \frac{q(1 + 2 \alpha)}{[2]_q([k]_q + \beta)}\geq \epsilon \big\}\bigg) = 0,$$
then
$$st-\lim_{n}\|{\mathcal{L}^{*(\alpha, \beta)}_n}(t;q_n,x)-x\|=0.$$
Similarly, by Lemma (\ref{c})(iii), we have
\begin{align*}
\Vert L_{n}^{\ast(\alpha,\beta)}(t^2;q_{n},x)-x^2\Vert &= \left\|\frac{[n]_q [m(n)]_q}{q([n]_q + {\beta})^2}x^2 + \frac{[n]_q \big[[3]_q+q\big((1+3\alpha)[2]_q+1 \big)\big]}{[3]_q([n]_q + {\beta})^2}x + \frac{q^2(1+ 3 \alpha + 3{\alpha}^2)}{[3]_q({[n]_q + {\beta}})^2}- x^2\right\|\\
&\leq \bigg|{ \frac{[n]_q [m(n)]_q}{q([n]_q + {\beta})^2}-1} \bigg|A^2+ \bigg|\frac{[n]_q [[3]_q+q((1+3\alpha)[2]_q+1)]}{[3]_q([n]_q + {\beta})^2}\bigg| A+ \bigg|\frac{q^2(1+ 3 \alpha + 3{\alpha}^2)}{[3]_q({[n]_q + {\beta}})^2}\bigg|
\end{align*}
\begin{align*}
~~~~\hspace{3.2cm}~~\leq {\mu}^2\left({\left(\frac{1}{q}-1\right)}+ \bigg(\frac{[n]_q(2+3\alpha)}{([n]_q + {\beta})^2}+\frac{[n]_q(q-(1+3\alpha))}{[3]_q([n]_q + {\beta})^2}\bigg) + \frac{q^2(1+ 3 \alpha + 3{\alpha}^2)}{[3]_q({[n]_q + {\beta}})^2}\right),
\end{align*}
where $\mu^2 = \max\{A^2,A,1\}=A^2.$\\
Now, if we choose
\begin{equation*}
\alpha _{n}= \left(\frac{1}{q}-1\right),
\end{equation*}%
\begin{equation*}
\beta _{n}=\frac{[n]_q(2+3\alpha)}{([n]_q + {\beta})^2}+\frac{[n]_q(q-(1+3\alpha))}{[3]_q([n]_q + {\beta})^2},
\end{equation*}%
\begin{equation*}
\gamma _{n}= \frac{q^2(1+ 3 \alpha + 3{\alpha}^2)}{[3]_q({[n]_q + {\beta}})^2},
\end{equation*}%
then by $(\ref{d}),$ we can write
\begin{equation} \label{2.1}
st-\lim_{n\rightarrow \infty }\alpha _{n}=0=st-\lim_{n\rightarrow \infty
}\beta _{n}=st-\lim_{n\rightarrow \infty }\gamma _{n}.
\end{equation}
Now for given $\epsilon >0$, we define the following four sets
\begin{equation*}
U=\{k: \Vert {\mathcal{L}^{*(\alpha, \beta)}_n}(t^{2};q_{k},x)-x^{2}\Vert \geq \epsilon \},
\end{equation*}%
\begin{equation*}
U_{1}=\{k:\alpha _{k}\geq \frac{\epsilon }{3{\mu}^2}\},
\end{equation*}%
\begin{equation*}
U_{2}=\{k:\beta _{k}\geq \frac{\epsilon }{3{\mu}^2}\},
\end{equation*}%
\begin{equation*}
U_{3}=\{k:\gamma _{k}\geq \frac{\epsilon }{3{\mu}^2}\}.
\end{equation*}%
It is obvious that $U\subseteq {U_{1}\cup U_{2}\cup U_{3}}$. Then, we obtain%
\begin{eqnarray*}
\delta \big(\{k &\leq & n:\Vert \mathcal{L}_{n}^{\ast(\alpha,\beta)}(t^{2};q_{n},x)-x^{2}\Vert \geq
\epsilon \}\big) \\&&~~~~~ \leq \delta \big(\{k\leq n:\alpha _{k}\geq \frac{\epsilon }{3{\mu}^2}\}\big) + \delta
\big(\{k\leq n:\beta _{k}\geq \frac{\epsilon }{3{\mu}^2}\}\big)+\delta \big(\{k\leq n:\gamma
_{k}\geq \frac{\epsilon }{3{\mu}^2}\}\big).
\end{eqnarray*}%
\newline
Using (\ref{2.1}), we get
\begin{equation*}
st-\lim_{n \rightarrow \infty }\Vert  {\mathcal{L}^{*(\alpha, \beta)}_n}(t^{2};q_{n},x)-x^{2}\Vert =0.
\end{equation*}
Since,
\begin{eqnarray*}
\Vert  {\mathcal{L}^{*(\alpha, \beta)}_n}(f;q_{n},x)-f\Vert \leq \Vert {\mathcal{L}^{*(\alpha, \beta)}_n}(t^{2};q_{n},x)-x^{2}\Vert +\Vert  {\mathcal{L}^{*(\alpha, \beta)}_n}
(t;q_{n},x)-x\Vert +\Vert {\mathcal{L}^{*(\alpha, \beta)}_n}(1;q_{n},x)-1\Vert ,\newline
\end{eqnarray*}
we get
\begin{eqnarray*}
st-\lim_{n\rightarrow \infty }\Vert {\mathcal{L}^{*(\alpha, \beta)}_n}(f;q_{n},x)-f\Vert
&\leq &  st-\lim_{n\rightarrow \infty }\Vert {\mathcal{L}^{*(\alpha, \beta)}_n}(t^{2};q_{n},x)-x^{2}\Vert\\
&&+~st-\lim_{n\rightarrow \infty }\Vert{\mathcal{L}^{*(\alpha, \beta)}_n}(t;q_{n},x)-x\Vert\\
&&+~st-\lim_{n\rightarrow \infty }\Vert {\mathcal{L}^{*(\alpha, \beta)}_n}(1;q_{n},x)-1\Vert ,
\end{eqnarray*}
which implies that
\begin{equation*}
st-\lim_{n\rightarrow \infty }\Vert  {\mathcal{L}^{*(\alpha, \beta)}_n}(f;q_{n},x)-f\Vert =0.
\end{equation*}
This completes the proof of theorem.

\section{Weighted statistical approximation}
Let $B_{x^2}[0,\infty)$ be set of all function $f$ defined on $[0,\infty)$ and satisfying the condition $ |f(x)|\leq M_{f}$ $ \rho(x)$, $M_f$ being a constant depending on $f$  and $\rho(x)=(1+x^2)\geq 1$ is called weighted function, it is continuous  on the positive real axis  and $\lim\limits_{x\rightarrow \infty} \rho(x)= \infty $. By $C_{x^2}[0,\infty)$, we denote the subspace of all continuous function belonging to $B_{x^2}[0,\infty)$. Also, $C_{x^2}^{*}[0,\infty)$ is subspace of all function $f\in C_{x^2}[0,\infty)$ for which $\displaystyle \lim_{x \rightarrow \infty} \frac{f(x)}{1+x^2}$ is finite. The norm on $C_{x^2}^{*}[0,\infty)$ is $\displaystyle \|f\|_{x^2}= \sup_{x\in[0,\infty)}\frac{|f(x)|}{1+x^2}.$
\begin{theorem}
Let $q=(q_n)$ be a sequence satisfying (\ref{d}) for $0< q_n \leq 1$. Then, for all non decreasing
functions $f\in C^{*}_{x^2}[0,\infty),$ we have
\begin{equation}
st-\lim_{n}\| {{\mathcal{L}^{*(\alpha, \beta)}_n}(f;q_{n},\cdot)-f\|}_{x^2}=0.
\end{equation}
\end{theorem}
\textbf{Proof.} By Lemma (\ref{c})(iii), we have  ${\mathcal{L}^{*(\alpha, \beta)}_n}(t^2; q_n, x) \leq Cx^2,$ where $C$ is a positive constant, ${\mathcal{L}^{*(\alpha, \beta)}_n}(f; q_n, x)$ is a sequence of positive linear operator acting from $C^{*}_{x^2}[0,\infty)$ to $B_{x^2}[0,\infty)$.\\
Using ${\mathcal{L}^{*(\alpha, \beta)}_n}(1; q_n, x) = 1,$ it is clear that\\
$$st-\lim_{n}\|{\mathcal{L}^{*(\alpha, \beta)}_n}(1;q_{n},x)-1\|_{x^2}=0.$$\\
Now, by Lemma (\ref{c})(ii), we have
$$\|{\mathcal{L}^{*(\alpha, \beta)}_n}(t; q_n, x)-x\|_{x^2}= \sup_{x\in [0,\infty)} \frac{|{\mathcal{L}^{*(\alpha, \beta)}_n}(t;q_n,x)-x|}{1+x^2} \leq \frac{\beta}{[n]_q + \beta}  + \frac{q(1 + 2 \alpha)}{[2]_q([n]_q + \beta)}.$$
Using (\ref{d}), we get
 $$st-\lim_{n}\left(\frac{\beta}{[n]_q + \beta}  + \frac{q(1 + 2 \alpha)}{[2]_q([n]_q + \beta)}\right)=0,$$
then $$st-\lim_{n}\|{\mathcal{L}^{*(\alpha, \beta)}_n}(t;q_n,x)-x\|_{x^2}=0.$$
Finally, by Lemma(\ref{c})(iii), we have 
\begin{align*}
\Vert {\mathcal{L}^{*(\alpha, \beta)}_n}(t^2;q_{n},x)-x^2\Vert_{x^2}
&\leq \left({ \frac{[n]_q [m(n)]_q}{([n]_q + {\beta})^2}-1} \right)\sup_{x\in[0,\infty)}\frac{x^2}{1+x^2}\\
&~~~+\left(\frac{[n]_q \big[[3]_q+q\big((1+3\alpha)[2]_q+1\big)\big]}{[3]_q([n]_q + {\beta})^2}\right)\sup_{x\in[0,\infty)}\frac{x}{1+x^2} + \left(\frac{q^2(1+ 3 \alpha + 3{\alpha}^2)}{[3]_q({[n]_q + {\beta}})^2}\right)\\
& \leq \left({ \frac{[n]_q [m(n)]_q}{([n]_q + {\beta})^2}-1} \right)+\left(\frac{[n]_q \big[[3]_q+q\big((1+3\alpha)[2]_q+1\big)\big]}{[3]_q([n]_q + {\beta})^2}\right)+ \left(\frac{q^2(1+ 3 \alpha + 3{\alpha}^2)}{[3]_q({[n]_q + {\beta}})^2}\right)\\
&=\left(\frac{1}{q_n}-1\right)+ \bigg(\frac{[n]_q(2+3\alpha)}{([n]_q + {\beta})^2}+\frac{[n]_q\big(q-(1+3\alpha)\big)}{[3]_q([n]_q + {\beta})^2}\bigg) + \left(\frac{q^2(1+ 3 \alpha + 3{\alpha}^2)}{[3]_q({[n]_q + {\beta}})^2}\right).
\end{align*}
If we choose
\begin{equation*}
\alpha _{n}= \left(\frac{1}{q_n}-1\right),
\end{equation*}%
\begin{equation*}
\beta _{n}=\frac{[n]_q(2+3\alpha)}{([n]_q + {\beta})^2}+\frac{[n]_q(q-(1+3\alpha))}{[3]_q([n]_q + {\beta})^2},
\end{equation*}%
\begin{equation*}
\gamma _{n}= \frac{q^2(1+ 3 \alpha + 3{\alpha}^2)}{[3]_q({[n]_q + {\beta}})^2},
\end{equation*}%
then by $(\ref{d}),$ we can write
\begin{equation}\label{3.1}
st-\lim_{n\rightarrow \infty }\alpha _{n}=0=st-\lim_{n\rightarrow \infty
}\beta _{n}=st-\lim_{n\rightarrow \infty }\gamma _{n}.
\end{equation}
Now for given $\epsilon >0$, we define the following four sets:
\begin{equation*}
S=\{k: \Vert {\mathcal{L}^{*(\alpha, \beta)}_n}(t^{2};q_{k},x)-x^{2}\Vert_{x^2} \geq \epsilon \},
\end{equation*}%
\begin{equation*}
S_{1}=\lbrace k:\alpha _{k}\geq \frac{\epsilon}{3}\rbrace,
\end{equation*}%
\begin{equation*}
S_{2}=\{k:\beta _{k}\geq \frac{\epsilon }{3}\},
\end{equation*}%
\begin{equation*}
S_{3}=\{k:\gamma _{k}\geq \frac{\epsilon }{3}\}.
\end{equation*}%
It is obvious that $S  \subseteq {S_{1}\cup S_{2}\cup S_{3}}$. Then, we obtain%
\begin{eqnarray*}
\delta \big(\{k &\leq &n:\Vert {\mathcal{L}^{*(\alpha, \beta)}_n}(t^{2};q_{n},x)-x^{2}\Vert \geq
\epsilon \}\big) \\
&&~~~~~~\leq \delta\big(\{k\leq n:\alpha _{k}\geq \frac{\epsilon }{3}\}\big)+\delta
\big(\{k \leq n:\beta _{k}\geq \frac{\epsilon }{3}\}\big)+\delta\big(\{k\leq n:\gamma
_{k}\geq \frac{\epsilon }{3}\}\big).
\end{eqnarray*}%
\newline
Using (\ref{3.1}), we get%
\begin{equation*}
st-\lim_{n\rightarrow \infty }\Vert {\mathcal{L}^{*(\alpha, \beta)}_n}(t^{2};q_{n},x)-x^{2}\Vert_{x^2} =0.
\end{equation*}
Since
\begin{eqnarray*}
&&\Vert{\mathcal{L}^{*(\alpha, \beta)}_n}(f;q_{n},x)-f\Vert_{x^2}\\
&&~~~~~~~~~\leq \Vert {\mathcal{L}^{*(\alpha, \beta)}_n}(t^{2};q_{n},x)-x^{2}\Vert_{x^2} +\Vert  {\mathcal{L}^{*(\alpha, \beta)}_n}
(t;q_{n},x)-x\Vert_{x^2}+\Vert {\mathcal{L}^{*(\alpha, \beta)}_n}(1;q_{n},x)-1\Vert_{x^2} ,\newline
\end{eqnarray*}%
we get
\begin{eqnarray*}
st-\lim_{n\rightarrow \infty }\Vert{\mathcal{L}^{*(\alpha, \beta)}_n}(f;q_{n},x)-f\Vert_{x^2} 
&\leq &st-\lim_{n\rightarrow \infty }\Vert{\mathcal{L}^{*(\alpha, \beta)}_n}(t^{2};q_{n},x)-x^{2}\Vert_{x^2}\\&& +~st-\lim_{n\rightarrow \infty }\Vert{\mathcal{L}^{*(\alpha, \beta)}_n}(t;q_{n},x)-x\Vert_{x^2}\\
&& + ~st-\lim_{n\rightarrow \infty }\Vert  {\mathcal{L}^{*(\alpha, \beta)}_n}(1;q_{n},x)-1\Vert_{x^2} ,
\end{eqnarray*}%
which implies that
\begin{equation*}
st-\lim_{n\rightarrow \infty }\Vert{\mathcal{L}^{*(\alpha, \beta)}_n}(f;q_{n},x)-f\Vert_{x^2} =0.
\end{equation*}
This completes the proof of the theorem.

\section{Rates of statistical convergence}
In this section, using the modulus of continuity, we study rates of statistical convergence of operator (\ref{b}) and Lipschitz functions are introduced.
\begin{lemma}\label{2.3}
Let $0<q<1$ and $a\in[0,bq],~ b>0.$ The inequality
\begin{equation}
\int_{a}^{b}|t-x|d_qt\leq \left( \int_{a}^{b}|t-x|^2d_qt\right)^{1/2} \left(\int_{a}^{b} d_qt \right)^{1/2}
\end{equation}
is satisfied.
\end{lemma}
~~ Let $C_B[0,\infty),$ the space of all bounded and continuous functions on $[0,\infty)$ and $x \geq 0.$
Then, for $\delta>0,$ the modulus of continuity of $f$ denoted by $\omega(f;\delta )$ is defined to be
$$\omega(f;\delta)= \sup_{|{t- x}|\leq {\delta}}|f(t)-f(x)|,~t\in[0,\infty).$$ 
It is known that $\lim\limits_{\delta\rightarrow 0}\omega(f ; \delta) = 0$ for $f\in C_B[0,\infty)$ and also, for any $\delta > 0$ and each $t,~x\geq 0,$ we have
\begin{equation}\label{2.2}
|f(t)-f(x)|\leq \omega(f;\delta)\left(1+\frac{|t-x|}{\delta}\right).
\end{equation}

\begin{theorem}
Let $(q_n)$ be a sequence satisfying (\ref{d}). For every non-decreasing $f\in C_B[0,\infty), ~x\geq 0$ and $n\in \mathbb{N},$ we have
$$|{\mathcal{L}^{*(\alpha, \beta)}_n}(f;q_{n},x)-f(x)| \leq 2\omega\big(f;\sqrt{\delta_n(x)} ~\big),$$
where
\begin{eqnarray}\label{2.4}
{\delta}_n(x)&=& \left (\frac{[n]_{q_n} [m(n)]_{q_n}}{{q_n}({[n]_{q_n} + \beta)}^2}+ 1- \frac{2[n]_{q_n}}{[n]_{q_n} + \beta} \right) x^2 \nonumber\\ &&+ \left(\frac{[n]_{q_n} \big[ [3]_{q_n} + {q_n}((1+3\alpha)[2]_{q_n}+1)\big] }{{ ([n]_{q_n} + \beta) }^2[3]_{q_n}} - \frac{2q_n(1+2\alpha)}{[2]_{q_n}([n]_{q_n} + \beta)} \right) x\nonumber  \\&&+ \left( \frac{{q_n}^2(1+ 3 \alpha+ 3 {\alpha}^2)}{[3]_{q_n} {([n]_{q_n} + \beta)}^2 } \right).
\end{eqnarray}
\end{theorem}
\textbf{Proof.}
Let non-decreasing $f\in C_B[0,\infty)$ and $x\geq 0$. Using linearity and positivity of the operators ${\mathcal{L}}^{*(\alpha, \beta)}_n$ and then applying (\ref{2.2}), we get for $\delta > 0$ and $n\in \mathbb{N}$ that
\begin{eqnarray*}
|{\mathcal{L}^{*(\alpha, \beta)}_n}(f;q_{n},x)-f(x)| &\leq &{\mathcal{L}^{*(\alpha, \beta)}_n}\big(|f(t)-f(x)|;q_{n},x\big)\\
&\leq & \omega(f,\delta)\big\{{\mathcal{L}^{*(\alpha, \beta)}_n}(1;q_{n},x)+ 
\frac{1}{\delta}{\mathcal{L}^{*(\alpha, \beta)}_n}(|t-x|;q_{n},x) \big\}.
\end{eqnarray*}
Taking ${\mathcal{L}^{*(\alpha, \beta)}_n}(1;q_{n},x) =1$ and then applying Lemma (\ref{2.3}) with $a=q\left(\frac{ [k]_q + q^{k-1}\alpha}{[n]_q+ \beta}\right)$ and $b={\frac{[k+1]_q + q^k \alpha}{[n]_q+ \beta}}$, we can write

\begin{align*}
|{\mathcal{L}^{*(\alpha, \beta)}_n}(f;q_{n},x)-f(x)|
&\leq  \omega(f;\delta)\bigg\{
1+\frac{1}{\delta}([n]_q + \beta) \sum_{k = 0}^{\infty} q^{k(k-1)/2} \frac{{D^k_q} (\phi_n (x))}{[k]_q ! } (-x)^k \\&~~~~~~~~~~~~\times{ \left(\int_{q \left( \frac{ [k]_q + q^{k-1}\alpha}{[n]_q+ \beta}\right)}^{\frac{[k+1]_q + q^k \alpha}{[n]_q+ \beta}}{ } {|q^{-k+1} t-x|}^2d_qt\right)}^{1/2}  {\left( \int_{q \left( \frac{ [k]_q + q^{k-1}\alpha}{[n]_q+ \beta}\right)}^{\frac{[k+1]_q + q^k \alpha}{[n]_q+ \beta}} d_qt\right)^{1/2}} \bigg\}.\\
~~~~~~~~~~~~~~~&\leq \omega(f;\delta)\bigg\{
1+\frac{1}{\delta} {\bigg(([n]_q + \beta)  \sum_{k = 0}^{\infty} q^{k(k-1)/2} \frac{{D^k_q} (\phi_n (x))}{[k]_q ! } (-x)^k  \int_{q \big( \frac{ [k]_q + q^{k-1}\alpha}{[n]_q+ \beta}\big)}^{\frac{[k+1]_q + q^k \alpha}{[n]_q+ \beta}}{ } {|q^{-k+1} t-x|}^2d_qt\bigg)}^{1/2} \\
&~~~~~~~~~~~
 \times  {\bigg(([n]_q + \beta) \sum_{k = 0}^{\infty} q^{k(k-1)/2} \frac{{D^k_q} (\phi_n (x))}{[k]_q ! } (-x)^k \int_{q \big( \frac{ [k]_q + q^{k-1}\alpha}{[n]_q+ \beta}\big)}^{\frac{[k+1]_q + q^k \alpha}{[n]_q+ \beta}} d_qt\bigg)^{1/2}}\bigg\} \\
~~~~~~~~~~~~~~~&\leq  \omega(f;\delta)\bigg\{1+\frac{1}{\delta}{\big({\mathcal{L}^{*(\alpha, \beta)}_n}((t-x)^2;q_{n},x)\big)}^{1/2}
{\big({\mathcal{L}^{*(\alpha, \beta)}_n}(1;q_{n},x)\big)}^{1/2}\bigg\}.
\end{align*}
Taking $q = (q_n),$ a sequence satisfying (\ref{d}), and using $\delta_n(x) = {\mathcal{L}^{*(\alpha, \beta)}_n}\big((t-x)^2;q_{n},x\big) $ and then
choosing $\delta = \delta_n(x)$ as in (\ref{2.4}), the theorem is proved.\\
\indent Observe that by the conditions in (\ref{d}), $st-\lim\limits_n \delta_n = 0.$ By (\ref{2.2}), we have
$$st-\lim\limits_n\omega(f;\delta_n) = 0.$$
This gives us the pointwise rate of statistical convergence of the operators  ${\mathcal{L}^{*(\alpha, \beta)}_n}(f;q_{n},x)$ to $f(x).$\\

\indent Now, we will study the rate of convergence of the operator ${\mathcal{L}^{*(\alpha, \beta)}_n}$ with the help of functions of the Lipschitz class $Lip_{M}(a)$, where $M>0$ and $0< a \leq 1$. Recall that a function
$f\in C_B[0,\infty)$ belongs to $Lip_{M}(a)$ if the inequality
\begin{equation}
|f(t)-f(x)| \leq M|t-x|^{a}; ~\forall ~ t, x\in[0,\infty)
\end{equation} holds.\\
Now, we have the following theorem.

\begin{theorem}
Let the sequence $q=(q_n)$ satisfy the condition given in (\ref{d}), and let $f\in Lip_{M}(a), ~x\geq 0$ with $0\leq a \leq 1$ and $M>0$. Then 
\begin{equation}
|{\mathcal{L}^{*(\alpha, \beta)}_n}(f;q_{n},x)-f(x)|\leq M \delta_n^{a/2}(x),
\end{equation}
where $\delta_n(x)$ is given as in (\ref{2.4}).
\end{theorem}
\textbf{Proof.}
Since ${\mathcal{L}^{*(\alpha, \beta)}_n}(f;q_{n},x)$ are linear positive operators and  $f\in Lip_{M}(a),$ on  $x\geq 0$ with $0<a < 1,$ we can write
\begin{align*}
|{\mathcal{L}^{*(\alpha, \beta)}_n}(f;q_{n},x)-f(x)|&\leq  {\mathcal{L}^{*(\alpha, \beta)}_n}(|f(t)-f(x)|;q_{n},x)\\
&\leq M {\mathcal{L}^{*(\alpha, \beta)}_n}(|t-x|^{a};q_{n},x).
\end{align*}
Now, we take $p=\frac{2}{a}$, $q=\frac{2}{2-a}$,  applying Lemma \ref{2.3} and H\"{o}lder's inequality, we obtain

\begin{align*}
|{\mathcal{L}^{*(\alpha, \beta)}_n}(f;q_{n},x)-f(x)|
&\leq  M\bigg\{
([n]_q + \beta) \sum_{k = 0}^{\infty} q^{k(k-1)/2} \frac{{D^k_q} (\phi_n (x))}{[k]_q ! } (-x)^k \\&~~~~~~\times{ \left(\int_{q \left( \frac{ [k]_q + q^{k-1}\alpha}{[n]_q+ \beta}\right)}^{\frac{[k+1]_q + q^k \alpha}{[n]_q+ \beta}}{ } {|q^{-k+1} t-x|}^2d_qt\right)}^{a/2}  {\left( \int_{q \left( \frac{ [k]_q + q^{k-1}\alpha}{[n]_q+ \beta}\right)}^{\frac{[k+1]_q + q^k \alpha}{[n]_q+ \beta}} d_qt\right)^{(2-a)/2}} \bigg\}\\
~~~~~~~~~~~~~~~&\leq M \bigg\{
 {\bigg(([n]_q + \beta)  \sum_{k = 0}^{\infty} q^{k(k-1)/2} \frac{{D^k_q} (\phi_n (x))}{[k]_q ! } (-x)^k  \int_{q \big( \frac{ [k]_q + q^{k-1}\alpha}{[n]_q+ \beta}\big)}^{\frac{[k+1]_q + q^k \alpha}{[n]_q+ \beta}}{|q^{-k+1} t-x|}^2d_qt\bigg)}^{a/2} \\
&~~~~~~
 \times  {\bigg(([n]_q + \beta) \sum_{k = 0}^{\infty} q^{k(k-1)/2} \frac{{D^k_q} (\phi_n (x))}{[k]_q ! } (-x)^k \int_{q \big( \frac{ [k]_q + q^{k-1}\alpha}{[n]_q+ \beta}\big)}^{\frac{[k+1]_q + q^k \alpha}{[n]_q+ \beta}} d_qt\bigg)^{(2-a)/2}}\bigg\} \\
~~~~~~~~~~~~~~~&\leq  M\bigg\{1+\frac{1}{\delta}{\big({\mathcal{L}^{*(\alpha, \beta)}_n}((t-x)^2;q_{n},x)\big)}^{1/2}
{\big({\mathcal{L}^{*(\alpha, \beta)}_n}(1;q_{n},x)\big)}^{1/2}\bigg\}.
\end{align*}
Taking $\delta_n(x)=\left({\mathcal{L}^{*(\alpha, \beta)}_n}\big((t-x)^2;q_{n},x\big)\right)$, as in (\ref{2.4}), we get
\begin{equation*}
|{\mathcal{L}^{*(\alpha, \beta)}_n}(f;q_{n},x)-f(x)|\leq M \delta_n^{a/2}(x).
\end{equation*}
Thus, the proof is complete.

\section{The construction of the bivariate operators of Kantorovich type}
The aim of this part is to construct the bivariate extension of the operator (\ref{b}), introduce the statistical convergence of the operators to the function $f$ and show the rate of statistical convergence of these operators. \\
 \indent $f:C\big([0,\infty)\times[0,\infty)\big)\rightarrow C\big([0,\infty)\times[0,\infty)\big)$
 and $0 < q_{n_1},\,q_{n_2}\leq 1$, let us define the bivariate case of operator (\ref{b}) as follows:

\begin{align} \label{2.5}
{\mathcal{L}^{*(\alpha, \beta)}_{{n_1},{n_2}}}(f;q_{n_1},q_{n_2},x,y)&=
([n_1]_{q_{n_1}} + \beta) ([n_2]_{q_{n_2}} + \beta)\nonumber\\ 
&~~~\times \sum_{k_1 = 0}^{\infty}\sum_{k_2 = 0}^{\infty} {q^{k_1(k_1-1)/2}_{n_1}} \frac{{D^{k_1}_{q_{n_1}}}(\phi_{n_1} (x))}{[k_1]_{q_{n_1}} ! } (-x)^{k_1} { q^{k_2(k_2-1)/2}_{n_2} \frac{{D^{k_2}_{q_{n_2}}} (\phi_{n_2} (x))}{[k_2]_{q_{n_2}} ! } (-x)^{k_2}} \nonumber\\
&~~~\times \displaystyle{\int_{q_{n_1} \big( \frac{ [k_1]_{q_{n_1}} + {q^{k_1-1}_{n_1}}\alpha}{[n_1]_{q_{n_1}}+ \beta}\big)}^{\frac{[k_1+1]_{q_{n_1}} + {q^{k_1}_{n_1}} \alpha}{[n_1]_{q_{n_1}}+ \beta}}}
 \int_{q_{n_2} \big( \frac{ [k_2]_{q_{n_2}} + {q^{k_2-1}_{n_2}}\alpha}{[n_2]_{q_{n_2}}+ \beta}\big)}^{\frac{[k_2+1]_{q_{n_2}} + {q^{k_2}_{n_2}} \alpha}{[n_2]_{q_{n_2}}+ \beta}}{ }{f\big({q^{-k_1+1}_{n_1}} t, ~ {q^{-k_2+1}_{n_2}} s\big)d_{q_{n_1}}t ~~d_{q_{n_2}}s}.
\end{align}

\indent In \cite{DU}, Erku\c{s} and Duman proved the statistical Korovkin type approximation theorem
for the bivariate linear positive operators to the functions in space $H_{{\omega}_2}.$\\
\indent Recently, Ersan and Do\u{g}ru \cite{ES1} obtained the statistical Korovkin type theorem and
lemma for the bivariate linear positive operators defined in the space $H_{{\omega}_2}$ as follows
 
\begin{theorem}\cite{ES1}
Let $\mathcal{D}_{n_1, n_2}$ be the sequence of linear positive operator acting from $H_{{\omega}_2} ({\mathbb{R}^2_+})$ into $C_B(\mathbb{R}_+),$ where $\mathbb{R}_+ = [0, \infty).$ Then, for any $f \in H_{{\omega}_2}$,
$$ st- \lim_{n_1, n_2} \| \mathcal{D}_{n_1, n_2}(f) - f \| = 0.$$
\end{theorem} 
 
\begin{lemma}
The bivariate operators defined in \cite{ES1} satisfy the following :
\begin{align*}
(i)~ \mathcal{D}_{n_1, n_2}(f_0; q_{n_1}, q_{n_2},x, y) &= q_{n_1}q_{n_2},\\
(ii)~ \mathcal{D}_{n_1, n_2}(f_1; q_{n_1}, q_{n_2},x, y) &= q_{n_1}q_{n_2} \frac{[n_1]_{q_{n_1}}}{[n_1+1]_{q_{n_1}}} \frac{x}{1+x},\\
(iii)~  \mathcal{D}_{n_1, n_2}(f_2; q_{n_1}, q_{n_2},x, y)& = q_{n_1}q_{n_2} \frac{[n_2]_{q_{n_2}}}{[n_2+1]_{q_{n_2}}} \frac{y}{1+y},\\
(iv)~  \mathcal{D}_{n_1, n_2}(f_3; q_{n_1}, q_{n_2},x, y)& = q^3_{n_1}q_{n_2} \frac{[n_1]_{q_{n_1}}[n_1-1]_{q_{n_1}}}{ {[n_1+1]^2}_{q_{n_1}} }  \frac{x^2}{(1+x)(1+q_{n_1}x)} +q_{n_1}q_{n_2} \frac{ [n_1]_{q_{n_1}}}{{[n_1+1]^2}_{q_{n_1}}} \frac{x}{(1+x)}\\&~~~~
 + q_{n_1}q^3_{n_2} \frac{[n_2]_{q_{n_2}}[n_2-1]_{q_{n_2}}}{{[n_2+1]^2}_{q_{n_2}}} \frac{y^2}{(1+y)(1+q_{n_2}y)} + q_{n_1}q_{n_2} \frac{[n_2]_{q_{n_2}}}{{[n_2+1]^2}_{q_{n_2}}} \frac{y}{1+y}.
\end{align*}
\end{lemma}

\indent ~~ In order to obtain the statistical convergence of  bivariate operator %
(\ref{2.5}), we need the following lemma.

\begin{lemma}\label{2.6} The bivariate operators defined in 
(\ref{2.5}) satisfy the following equalities:
\begin{align*}
{\mathcal{L}_{n_1, n_2}^{*(\alpha, \beta)}}(e_0;q_{n_1},q_{n_2},x,y) &= 1,\\
{\mathcal{L}_{n_1, n_2}^{*(\alpha, \beta)}}(e_1;q_{n_1},q_{n_2},x, y)&= \frac{[n_1]_{q_{n_1}}}{[n_1]_{q_{n_1}} + \beta} x + \frac{q_{n_1}(1 + 2 \alpha)}{[2]_{q_{n_1}}([n_1]_{q_{n_1}} + \beta)},\\
{\mathcal{L}_{n_1, n_2}^{*(\alpha, \beta)}}(e_2;q_{n_1},q_{n_2}, x, y)&= \frac{[n_2]_{q_{n_2}}}{[n_2]_{q_{n_2}} + \beta} y + \frac{q_{n_2}(1 + 2 \alpha)}{[2]_{q_{n_2}}([n_2]_{q_{n_2}} + \beta)},\\
{\mathcal{L}_{n_1, n_2}^{*(\alpha, \beta)}}(e_3;q_{n_1},q_{n_2}, x, y)&= \frac{[n]_{q_{n_1}} [m(n)]_{q_{n_1}}}{{q_{n_1}}([n]_{q_{n_1}} + {\beta})^2}x^2 + \frac{[n]_{q_{n_1}} \big[[3]_{q_{n_1}}+{q_{n_1}}((1+3\alpha)[2]_{q_{n_1}}+1)\big]}{[3]_{q_{n_1}}([n]_{q_{n_1}} + {\beta})^2}x + \frac{{q}^2_{n_1}(1+ 3 \alpha + 3{\alpha}^2)}{[3]_{q_{n_1}}({[n]_{q_{n_1}} + {\beta}})^2}\\ &~~~+ \frac{[n]_{q_{n_2}} [m(n)]_{q_{n_2}}}{([n]_{q_{n_2}} + {\beta})^2}y^2 + \frac{[n]_{q_{n_2}} \big[[3]_{q_{n_2}}+{q_{n_2}}((1+3\alpha)[2]_{q_{n_2}}+1)\big]}{[3]_{q_{n_2}}([n]_{q_{n_2}} + {\beta})^2}y + \frac{{q}^2_{n_2}(1+ 3 \alpha + 3{\alpha}^2)}{[3]_{q_{n_2}}({[n]_{q_{n_2}} + {\beta}})^2}.
\end{align*}
\end{lemma}
\textbf{Proof.} The proof can be obtained similar to the proof of bivariate operator in \cite{ES1}. So, we shall omit this proof.\\

\indent ~~ Let $q = (q_{n_1})$ and  $q = (q_{n_2})$ be the sequence that converges statistically to $1$ but does not
converge in ordinary sense, so for $ 0< q_{n_1},\,q_{n_2} \leq 1,$ it can be written as
\begin{equation}\label{2.7}
 st-\lim_{n_1} q_{n_1}= st-\lim_{n_2} q_{n_2}=1.
\end{equation}
Now, under the condition in (\ref{2.7}), let us show the statistical convergence of bivariate operator
(\ref{2.5}) with the help of the proof of Theorem \ref{A1}.

\begin{theorem}
Let $q= (q_{n_1})$ and $q=(q_{n_2})$ be sequence satisfying (\ref{2.7}) for $0<q_{n_1}, q_{n_2}\leq  1,$ and let 
${\mathcal{L}^{*(\alpha, \beta)}_{n_1,n_2}}$ be sequence of  linear positive operator from $C(K)$ into $C(K)$ given by (\ref{2.5}). Then, for any function $f\in C(K_1\times K_1)\subset C(K\times K)$ and $x\in K_1\times K_1\subset K\times K,$ where  $K=[0,\infty)\times [0,\infty)$, $K_1=[0, A]\times [0, A],$ we have 
\begin{equation}
st-\lim_{n_1,n_2}\parallel{\mathcal{L}^{*(\alpha, \beta)}_{n_1,n_2}}(f)-f \parallel_{C(K_1\times K_1)}=0.
\end{equation}
\end{theorem}
\textbf{Proof.} Using Lemma (\ref{2.6}), the proof can be obtained similar to the proof of Theorem \ref{A1}. So, we shall omit this proof.

\section{Rates of convergence of bivariate operators}
Let $K=[0,\infty)\times [0,\infty).$ Then  the sup norm on $C_B(K)$ is given by 
$$ \|f\|= \sup_{(x,y)\in K}|f(x,y)|,~~~ f\in C_B(K).$$
We consider the modulus of continuity $\omega(f;\delta_1,\delta_2)$, where $\delta_1,~\delta_2 >0,$ for bivariate case given by 
\begin{equation}
\omega(f;\delta_1,\delta_2)=\{\sup |f(x',y')-f(x,y)|:(x',y'),\,(x,y)\in K ~\text {and}~ |x'-x|\leq \delta_1, ~ |y'-y|\leq\delta_2\}.
\end{equation}
It is clear that a necessary and sufficient condition for a function $f\in C_B(K)$ is
$$\lim_{\delta_1,\,\delta_2\rightarrow 0}\omega(f;\delta_1,\delta_2)=0$$
and $\omega(f;\delta_1,\delta_2)$ satisfy the following condition:
\begin{equation}\label{2.8}
|f(x',y')-f(x,y)|\leq \omega(f;\delta_1,\delta_2)\left(1+\frac{|x'-x|}{\delta_1}\right)\left(1+\frac{|y'-y|}{\delta_2}\right)
\end{equation}
for each $f\in C_B(K)$. Then observe that any function in $C_B(K)$ is continuous and bounded
on $K$. Details of the modulus of continuity for bivariate case can be found in \cite{GA}.\\
\indent Now, the rate of statistical convergence of bivariate operator (\ref{2.5}) by means of modulus
of continuity in $f \in C_B(K)$ will be given in the following theorem.

\begin{theorem}
Let $q= (q_{n_1})$ and $q=(q_{n_2})$ be sequence satisfying (\ref{2.7}). So, we have
\begin{equation}
\vert{\mathcal{L}}^{*(\alpha, \beta)}_{n_1,n_2}(f;q_{n_1},q_{n_2},x,y)-f(x,y)\vert\leq 4\omega\big(f;\sqrt{\delta_{n_1}(x)},\sqrt{\delta_{n_2}(y)}~\big),
\end{equation}
where
\begin{align}\label{2.9}
{\delta}_{n_1}(x)&= \left (\frac{[{n_1}]_{q_{n_1}} [m({n_1})]_{q_{n_1}}}{{q_{n_1}}({[{n_1}]_{q_{n_1}} + \beta)}^2}+ 1- \frac{2[{n_1}]_{{q_{n_1}}}}{([{n_1}]_{{q_{n_1}}} + \beta)} \right) x^2 \nonumber\\
&~~~~~~~+ \left(\frac{[{n_1}]_{q_{n_1}} \big[ [3]_{q_{n_1}}+ {q_{n_1}}((1+3\alpha)[2]_{q_{n_1}}+1)\big]}{{ ([{n_1}]_{q_{n_1}} + \beta) }^2[3]_{q_{n_1}}} - \frac{2{q_{n_1}}(1+2\alpha)}{[2]_{q_{n_1}}([{n_1}]_{q_{n_1}} + \beta)} \right) x\nonumber  \\
&~~~~~~~+ \left( \frac{{q}^2_{n_1}(1+ 3 \alpha+ 3 {\alpha}^2)}{[3]_{q_{n_1}} {([{n_1}]_{q_{n_1}} + \beta)}^2 } \right), 
\end{align}
\begin{align}\label{2.10}
{\delta}_{n_2}(y)&= \left (\frac{[n_2]_{q_{n_2}} [m(n_2)]_{q_{n_2}}}{{q_{n_2}}({[n_2]_{q_{n_2}} + \beta)}^2}+ 1- \frac{2[n_2]_{q_{n_2}}}{([n_2]_{q_{n_2}} + \beta)} \right) y^2 \nonumber\\
&~~~~~~~+ \left(\frac{[n_2]_{q_{n_2}} \big[ [3]_{q_{n_2}}+ {q_{n_2}}((1+3\alpha)[2]_{q_{n_2}}+1)\big] }{{ ([n_2]_{q_{n_2}} + \beta) }^2[3]_{q_{n_2}}} - \frac{2{q_{n_2}}(1+2\alpha)}{[2]_{q_{n_2}}([n_2]_{q_{n_2}} + \beta)} \right) y \nonumber  \\
&~~~~~~~+  \left( \frac{{q}^2_{n_2}(1+ 3 \alpha+ 3 {\alpha}^2)}{[3]_{q_{n_2}} {([n_2]_{q_{n_2}} + \beta)}^2 } \right). 
\end{align}
\end{theorem}
\textbf{Proof.} By using the condition in (\ref{2.8}), for  $\delta_{n_1},\delta_{n_2} > 0 $
 and $n\in \mathbb{N}$, we get 
\begin{align*}
&\vert{\mathcal{L}}^{*(\alpha, \beta)}_{n_1,n_2}(f;q_{n_1},q_{n_2},x,y)-f(x,y)\vert \\
&~~~~~~~\leq {\mathcal{L}}^{*(\alpha, \beta)}_{n_1,n_2}(\vert f(x',y')-f(x,y)\vert; q_{n_1},q_{n_2},x,y)\\
&~~~~~~~\leq \omega(f;\delta_{n_1}(x),\delta_{n_2}(y)) \bigg({\mathcal{L}}^{*(\alpha, \beta)}_{n_1,n_2}(f_0; q_{n_1},q_{n_2},x,y)+\frac{1}{\delta_{n_1}}{\mathcal{L}}^{*(\alpha, \beta)}_{n_1,n_2}(\vert x'- x\vert; q_{n_1},q_{n_2},x,y)\bigg) \\
&~~~~~~~~~~~\times \bigg({\mathcal{L}}^{*(\alpha, \beta)}_{n_1,n_2}(f_0; q_{n_1},q_{n_2},x,y)+\frac{1}{\delta_{n_2}}{\mathcal{L}}^{*(\alpha, \beta)}_{n_1,n_2}(\vert y'- y\vert; q_{n_1},q_{n_2},x,y)\bigg)
\end{align*}
If the Cauchy-Schwarz inequality is applied, we have
$${\mathcal{L}}^{*(\alpha, \beta)}_{n_1,n_2}(\vert x'- x\vert; q_{n_1},q_{n_2},x,y)\leq \bigg({\mathcal{L}}^{*(\alpha, \beta)}_{n_1,n_2}((x'- x)^2; q_{n_1},q_{n_2},x,y)\bigg)^{1/2}\bigg({\mathcal{L}}^{*(\alpha, \beta)}_{n_1,n_2}(f_0; q_{n_1},q_{n_2},x,y)\bigg)^{1/2}. $$
So, if it is substituted in the above equation, the proof is completed.\\

\indent At last, the following theorem represents the rate of statistical convergence of bivariate
operator (\ref{2.5}) by means of Lipschitz $Lip_M(\alpha_1,\alpha_2)$ functions for the bivariate case, where
$f\in C_B[0,\infty)$ and $M >0$ and $0 < \alpha_1\leq 1$, $0 < \alpha_2 \leq 1$, then let us define $Lip_M(\alpha_1,\, \alpha_2)$ as

$$\vert f(x',y')-f(x,y)\vert \leq M {\vert x'- x\vert}^{\alpha_1}{\vert y'- y\vert}^{\alpha_2};~~~~~~ \forall \,  x,x',y,y'\in[0,\infty).$$
We have the following theorem.

\begin{theorem}
Let $q= (q_{n_1})$ and $q=(q_{n_2})$ be sequence satisfying the condition given in (\ref{2.7}), and let  $Lip_M(\alpha_1,\alpha_2),$ $x\geq 0$ and $0 <\alpha_1 \leq 1,~ 0<\alpha_2\leq 1$. Then

\begin{equation}
\vert{\mathcal{L}}^{*(\alpha, \beta)}_{n_1,n_2}(f;q_{n_1},q_{n_2},x,y)-f(x,y)\vert\leq M \, {\delta_{n_1}^{\alpha_1/2}(x)}\, {\delta_{n_2}^{\alpha_2/2}(y)},
\end{equation}
where $\delta_{n_1}(x)$ and  $\delta_{n_2}(y)$ are defined in (\ref{2.9}), (\ref{2.10}).
\end{theorem}  
\textbf{Proof.} Since ${\mathcal{L}}^{*(\alpha, \beta)}_{n_1,n_2}(f;q_{n_1},q_{n_2},x,y)$  are linear positive operators and $ f\in Lip_M(\alpha_1,\alpha_2),~  x\geq 0$ and $0 <\alpha_1 \leq 1, 0<\alpha_2\leq 1,$ we can write
\begin{align*}
\vert{\mathcal{L}}^{*(\alpha, \beta)}_{n_1,n_2}(f;q_{n_1},q_{n_2},x,y)-f(x,y)\vert
&\leq {\mathcal{L}}^{*(\alpha, \beta)}_{n_1,n_2}(\vert f(x',y')-f(x,y)\vert; q_{n_1},q_{n_2},x,y)\\
&\leq  M {\mathcal{L}}^{*(\alpha, \beta)}_{n_1,n_2}(|x'-x|^{\alpha_1} |y'-y|^{\alpha_2}; q_{n_1},q_{n_2},x,y)\\
& = M  {\mathcal{L}}^{*(\alpha, \beta)}_{n_1,n_2}(|x'-x|^{\alpha_1} ; q_{n_1},q_{n_2},x,y)\, {\mathcal{L}}^{*(\alpha, \beta)}_{n_1,n_2}(|y'-y|^{\alpha_2}; q_{n_1},q_{n_2},x,y).
\end{align*} 
 If we take $p_1 = \frac{2}{\alpha_1},\, p_2 = \frac{2}{\alpha_2}, \, q_1 = \frac{2}{2-\alpha_1},\,  q_2 = \frac{2}{2-\alpha_2},$ applying H\"{o}lder's inequality, we obtain\\
\begin{align*}
\vert{\mathcal{L}}^{*(\alpha, \beta)}_{n_1,n_2}(f;q_{n_1},q_{n_2},x,y)-f(x,y)\vert
&\leq \bigg({\mathcal{L}}^{*(\alpha, \beta)}_{n_1,n_2}\big((x'- x)^{\alpha_1}; q_{n_1},q_{n_2},x,y\big)\bigg)^{\alpha_1/2} \bigg({\mathcal{L}}^{*(\alpha, \beta)}_{n_1,n_2}(f_0; q_{n_1},q_{n_2},x,y)\bigg)^{(2-\alpha_1)/2}\\
& ~~~\times \bigg({\mathcal{L}}^{*(\alpha, \beta)}_{n_1,n_2}\big((y'- y)^{\alpha_2}; q_{n_1},q_{n_2},x,y\big)\bigg)^{\alpha_2/2}
\bigg({\mathcal{L}}^{*(\alpha, \beta)}_{n_1,n_2}(f_0; q_{n_1},q_{n_2},x,y)\bigg)^{(2-\alpha_2)/2}\\
&= M {\delta^{\alpha_1/2}_{n_1}}(x){\delta^{\alpha_2/2}_{n_2}}(y).
\end{align*}
Which is the required result.\\

\hspace{-.9cm}
\textbf{Conclusion}\\
Our proposed family of integral operators ${\mathcal{L}^{*(\alpha, \beta)}_n}$ are generalization of summation-integral type operators. The results established here are more general rather than the results of any other previous proved lemmas and theorems. The strong convergence in weighted spaces is highlighted and Bivariate generalization also established for said operators. Some special cases are also considered. Problems considered in this paper may open further research opportunities in these fields. The researchers and professionals working or intend to work in the areas of analysis and its applications will find this research article to be quite useful. 

\vspace{0.5cm}
\hspace{-.9cm}
\textbf{Conflict of Interests}\\
The authors declare that there is no conflict of interests regarding the publication of this paper.

\vspace{0.5cm}
\hspace{-.9cm}
\textbf{Acknowledgement}\\
The authors would like to express their deep gratitude to the anonymous learned referee(s) and the editor for their valuable suggestions and constructive comments, which resulted in the subsequent improvement of this research article. Special thanks are due to our great Master and friend academician Prof. Hari Mohan Srivastava, Editor of Filomat for kind cooperation, smooth behaviour during communication and for his efforts to send the reports of the manuscript timely. The authors are also grateful to all the editorial board members and reviewers of esteemed journal i.e. Filomat. The second author PS is thankful to Department of Applied Mathematics and humanities, SVNIT, Surat (Gujarat) to carry out her research work (Ph.D. in Full-time Institute Research (FIR) category) under the supervision of Dr. Vishnu Narayan Mishra. The first author VNM acknowledges that this project was supported by the Cumulative Professional Development Allowance (CPDA), SVNIT, Surat (Gujarat), India. All the authors carried out the proof of Lemmas and Theorems. Each author contributed equally in the development of the manuscript. VNM conceived of the study and participated in its design and coordination. All the authors read and approved the final version of manuscript.


\end{document}